\documentclass[12pt]{article}
\newtheorem{theorem}{Theorem}
\newtheorem{definition}{Definition}
\newtheorem{lemma}{Lemma}
\newtheorem{corollary}{Corollary}
\usepackage{amssymb}
\begin{document}

\title {Cycles of weight divisible by $k$}
\author{Ajit A. Diwan \\
Department of Computer Science and Engineering,\\
Indian Institute of Technology Bombay,\\
Mumbai 400076.\\
email: aad@cse.iitb.ac.in }

\maketitle

\abstract{ A weighted (directed) graph is a (directed) graph with integer weights assigned to its 
vertices and edges. The weight of a subgraph is the sum of weights of vertices and edges in the 
subgraph. The problem of determining the largest order $f(k)$ of a weighted complete directed 
graph that does not contain a directed cycle of weight divisible by $k$, for an integer $k \ge 2$,
was raised by Alon and Krivelevich [J. Graph Theory 98 (2021) 623--629]. They showed that $f(k)$ is 
$O(k\log k)$ and $f(k) \le 2k-2$ if $k$ is prime. The best bounds known to us 
are $f(k) \le 2k-2$ for all $k$ and $f(k) < (3k-1)/2$ for prime $k$. It is also known that
$f(k) \ge k$ and this is believed to be the correct value. We prove that $f(k) < k+2\Omega(k)$, 
where $\Omega(k)$ is the number of prime factors, not necessarily distinct, in the prime factorization 
of $k$.

We also show that any weighted undirected graph of minimum degree $2k-1$ contains a cycle
of weight divisible by $k$. This result is proved in the more general setting in which the weights
are from a finite abelian group of order $k$, and the cycle has weight equal to the group
identity. We conjecture that this holds for undirected graphs with minimum degree $k+1$.}

\section{Introduction}

We consider only simple graphs that may be directed or undirected. A weight function on a
graph $G$ is a function $w : V(G) \cup E(G) \rightarrow \mathbb{Z}$ that assigns to each
vertex and edge an integer weight. We call $G$ together with a weight function $w$ a
$w$-weighted graph and omit the $w$ if it is understood from the context. The weight of a subgraph 
of a weighted graph is the sum of weights of the vertices and edges in the subgraph.

The central question in zero-sum Ramsey theory is to find the smallest order of a weighted
complete (directed) graph that ensures that it contains a subgraph of a particular kind whose 
weight is divisible by $k$, for a given integer $k \ge 2$. This is equivalent to considering
the weights to be in the additive group of integers modulo $k$ and requiring the
weight of the subgraph to be 0. We call such a subgraph a zero subgraph. A somewhat old survey 
of results of this kind can be found in~\cite{car}.

The simplest possible subgraph to consider is a (directed) cycle. The question of
determining the largest order $f(k)$ of a weighted complete directed graph that does not
contain a cycle of weight divisible by $k$ was raised by Alon and Krivelevich~\cite{AK}.
They proved that $f(k)$ is $O(k\log k)$ and $f(k) \le 2k-2$ for a prime $k$. This was
improved in~\cite{MS} where it is shown that $f(k) < 8k$ for all $k$ and $f(k) < (3k-1)/2$ for prime $k$. 
They also showed the easy lower bound $f(k) \ge k$ and asked the question whether it is the
optimal bound. In~\cite{LM}, where the same problem is considered for the group $\mathbb{Z}_p^k$, 
it is mentioned that in a forthcoming work, the optimal bound $k$ has been proved. However,
we have been unable to access the manuscript. The best upper bound that we know for general $k$ is 
$f(k) \le 2k-2$ for all $k$. This is proved in~\cite{AAG,BBK}, and also holds in the more general 
setting of group weighted graphs, even for nonabelian groups. Here we prove that $f(k) < k + 2\Omega(k)$,
for $k \ge 2$, where $\Omega(k)$ is the number of prime factors, not necessarily distinct, in 
the prime factorization of $k$. This is within one of the optimal bound for prime $k$, and
is asymptotically equal to the optimal for large $k$.

We also consider the case of undirected graphs. In the directed case, when considering directed
cycles, weights on vertices can be ignored, since the weight of a vertex can be added to
the weights of all edges directed out of the vertex, without changing the weight of any cycle.
However, this is not true if the graph is undirected. If the weight of a vertex is even, then
half the weight can be added to the weights of all edges incident with the vertex, but such
a reduction is not possible if the weight is odd. Also, if an undirected edge has weight 0,
replacing it by two oppositely directed edges of weight 0 creates a zero cycle in the directed
graph, but it is not a cycle in the undirected graph. However, this can be taken care of
by requiring the directed cycle to have length at least 3.  With this restriction, we show that
the bound on the order increases by just one.

A major difference in the undirected case though is that the result holds not just for 
complete graphs but also for graphs with large enough minimum degree. We show that
every weighted undirected graph of minimum degree at least $2k-1$, for $k \ge 2$, contains
a cycle of weight divisible by $k$. We conjecture that the correct bound is in fact $k+1$
and give examples of weighted undirected graphs with minimum degree $k$ that do not
contain a cycle of weight divisible by $k$. In particular, there is a weighted complete graph
of order $k+1$ that does not contain a cycle of weight divisible by $k$. The proof of the upper 
bound holds in the more general setting of weights from an abelian group of order $k$.

\section{Directed graphs}

In this section, we consider weighted complete directed graphs. We assume the weights are from the
set $\mathbb{Z}_k$ of integers modulo $k$, for some $k \ge 2$. It is sufficient to consider the case
when all vertex weights are 0, since as mentioned earlier, the weight of a vertex can be added to the weights 
of all edges directed out of the vertex, without changing the weight of any cycle. We will therefore
assume the weight of every vertex is 0. If $d$ is any proper divisor of
$k$, and $0 \le i < d$,  the congruence class $i$ modulo $d$ is the set of all integers in
$\mathbb{Z}_k$ that are congruent to $i$ modulo $d$. A congruence class modulo $d$ is simply
a congruence class $i$ modulo $d$ for some $0 \le i < d$. If $A \subset \mathbb{Z}_k$ and
$a \in \mathbb{Z}_k$,  we denote by $A+a$ the set of all elements in $\mathbb{Z}_k$ that
can be obtained by adding $a$ to an element in $A$. Note that $|A+a| = |A|$. Let $\Omega(k)$ denote
the number of prime factors, not necessarily distinct, in the prime factorization of $k$.
If $d > 1$ is any proper divisor of $k$, then $\Omega(k) = \Omega(k/d) + \Omega(d)$, hence
$\Omega(k) \ge \Omega(k/d)+1$.

A set $A \subset \mathbb{Z}_k$ is said to be a near arithmetic progression if $2 \le |A| \le k-2$
and there exists a subset $B \subset A$ with $|B|= |A|-1$ and some nonzero element $a \in
\mathbb{Z}_k$ such that $B+a \subset A$. The main property of $\mathbb{Z}_k$ that we use is given 
by the following lemma. 

\begin{lemma}
\label{inc}
Let $A \subset \mathbb{Z}_k$ be a near arithmetic progression. Then at least
one of the following statements holds for $A$.
\begin{enumerate}
\item
There exists a proper divisor $d > 1$ of $k$, such that for any subset $A' \subset A$ with $|A'| =
|A|-1$ and any $x \in \mathbb{Z}_k$, if $A'+x \subset A$ then $x$ is a multiple of $d$.
\item
There exists an element $a \in \mathbb{Z}_k$ such that $gcd(a,k) = 1$ and for any subset $A' \subset A$ with 
$|A'| = |A|-1$ and any $x \in \mathbb{Z}_k$, if $A'+x \subset A$ then $x \in \{0,a,-a\}$.
\end{enumerate}
\end{lemma}

\noindent {\bf Proof:}  Since $A$ is a near arithmetic progression, by definition, there exists a
subset $B \subset A$ and a nonzero element $a \in \mathbb{Z}_k$ such that $|B| = |A|-1$ and 
$B+a \subset A$. Let $b$ be the element in $A \setminus B$. 

\noindent {\bf Case 1.} Suppose $b \not\in B+a$.

This implies $B = B+a$ so for every element $x \in B$, $x +ia$ is in $B$ for all $i  \ge 0$,
hence $x+gcd(a,k) \in B$. This implies $B$ is the union of some congruence classes modulo $gcd(a,k)$.
Let $d$ be the smallest divisor of $k$ such that $B$ is the union of congruence classes modulo $d$.
We show that in this case $A$ satisfies the first property. Since $|A| < k$, we have $d > 1$. Also, since 
$|A| > 1$, $B$ contains at least one congruence class modulo $d$. 

Let $A' \subset A$ with $|A'| = |A|-1$ and suppose $A'+x \subset A$ for some $x$. Suppose $x$ is not 
a multiple of $d$ and let $x = c$ modulo $d$. First suppose $A'= B$. Let $S = B \cap \mathbb{Z}_d$.
If for some element $i \in S$, $i+c$ modulo $d$ is not in $S$, the congruence class $i+c$ modulo
$d$ is contained in $B+x$ but not in $B$. Hence $A$ can contain at most one element of that class, a 
contradiction. Therefore for all $i \in S$, $i+c$ modulo $d$ is also in $S$. However, this implies $B$ is the
union of congruence classes modulo $gcd(c,d)$, contradicting the choice of $d$. Therefore in
this case $x$ must be a multiple of $d$. 

Suppose $A' \neq B$ which implies $b \in A'$ and hence $A'+x$ contains exactly one number congruent to 
$b+c$ modulo $d$. However, for any congruence class modulo $d$ that does not contain $b$, either $A$ 
contains all elements in the class or none. If $d \le k/3$, each congruence class contains at least 
3 elements, which implies there are at least 2 elements in $A$ that are not in $A'+x$, contradicting the
assumption that $|A'| = |A|-1$ and $A'+x \subset A$. This implies that if $d \le k/3$, $x$ must be a multiple 
of $d$ in this case also. Thus if $d \le k/3$, $A$ satisfies the first property with $d = gcd(a,k)$.

The remaining case is if $d = k/2$ and $A'$ contains $b$. Let $A'' = A \cup \{b+k/2\}$. Then $A''$ is a 
union of congruence classes modulo $k/2$. Let $S = A'' \cap \mathbb{Z}_{k/2}$. Suppose there exists an element 
$i \in S$ such that $i+c$ modulo $k/2$ is not in $S$. This implies $A'$ cannot contain any of the elements 
$i, i+k/2$, otherwise $A'+x$ contains an element not in $A$. However, since $b \in A'$, $b+k/2 \not\in A'$ 
and $|A'| = |A''|-2$, this gives a contradiction. Therefore for every element $i \in S$, $i+c$ modulo $k/2$ 
is in $S$, and hence $S$ is the union of congruence classes modulo $gcd(c,k/2)$. This implies $A''$ is also
a union of congruence classes modulo $gcd(c,k/2)$. Let $d'$ be the smallest divisor of $gcd(c,k/2)$ such 
that $A''$ is a union of congruence classes modulo $d'$. Note that since $c < k/2$, $d' \le k/4$ and
each congruence class modulo $d'$ contains at least 4 elements. Also, since $|A''| < k$, $d' > 1$. We claim 
that $x$ must be a multiple of $d'$. Again, if not, and $x = c'$ modulo $d'$, for every congruence class 
$i$ modulo $d'$ in $A''$, the congruence class $i+c'$ must be in $A''$, otherwise we get at least 4 elements 
in $A''$ that are not in $A'$, a contradiction. On the other hand, if this holds, it contradicts the choice of 
$d'$. Therefore in this case $x$ must be multiple of $d'$. Since $d'$ is a divisor of $d$, in all cases
$x$ must be a multiple of $d'$.  This implies $A$ satisfies the first property. An example of such a case
is if $k = 8$, $A = \{0,2,4\}$, $B = \{0,4\}$ and $a = 4$. Here $d=4$ but $d' = 2$ and for the set
$A' = \{0,2\}$, $A'+2 \subset A$.

\noindent {\bf Case 2.}  Suppose $b \in B+a$.

In this case, there must be exactly one element $b' \in A$ such that $b' \not\in B+a$ and $b' \neq b$,
hence $b' \in B$. Consider the longest arithmetic progression contained in $A$ that starts with $b'$ and has 
common difference $a$. The last element in this sequence must be $b$, otherwise if it it an element in $B$,
we can add $a$ to it to get a longer progression.  Note that the element added cannot be $b'$ since
$b' \not\in B+a$. Let $A_2$ be the elements of $A$ in this arithmetic progression and let $A_1 = A \setminus
A_2$. Then $A_1+a = A_1$ and $A_1$ is the union (possibly empty) of congruence classes modulo $d = gcd(a,k)$.
The arithmetic progression $A_2$ is contained in one congruence class modulo $d$.  

Suppose $d=1$, which implies $A_1$ is empty. Since $|A| \le k-2$, the elements $b'-a$, $b'-2a$ 
are not in $A$. In this case we show that $A$ satisfies the second property. Let $A'$ be a subset
of $A$ with $|A'| = |A|-1$ such that $A'+x \subset A$. Let the elements in $A_2$ be $b'+ia$ for
$0 \le i < l$, where $l = |A|$. At least one of the elements $b'$ or $b'+a$ must be contained in $A'+x$,
otherwise there are 2 elements in $A$ not in $A'+x$, contradicting the assumption that $|A'| = |A|-1$
and $A'+x \subset A$. Suppose $b' \in A'+x$. Then for some $0 \le i < l$ we have $b' = b'+ia + x$.
If $i = 0$, then $x = 0$ and if $i = 1$ then $x = -a$. Suppose $i \ge 2$. Then if $b'+(i-1)a \in A'$,
then $b'+(i-1)a + x = b'-a \not\in A$, contradicting the assumption that $A'+x \subset A$. If
$b'+(i-1)a \not\in A'$ then $b'+(i-2)a \in A'$ and $b'+(i-2)a+x = b'-2a \not\in A$, again a contradiction.
Suppose $b' \not\in A'+x$ but $b'+a \in A'+x$. Then for some $0 \le i < l$, $b'+a = b'+ia+x$. If $i = 0$,
then $x = a$ and if $i = 1$ then $x = 0$. If $i \ge 2$, then since
$b', b'-a \not\in A'+x$, neither $b'+(i-1)a$ nor $b'+(i-2)a$ are in $A'$, contradicting the assumption
that $|A'| = |A|-1$. Therefore, for any such subset $A'$ and element $x$, we must have $x \in \{0,a,-a\}$
with $gcd(a,k) = 1$.

Suppose $d > 1$. Let $A''$ be the union of the congruence classes $i$ modulo $d$, for all $i \in \mathbb{Z}_d$ 
such that $A$ contains an element congruent to $i$ modulo $d$. Suppose $|A''| < k$ and let $d'$ be the
smallest divisor of $d$ such that $A''$ is the union of congruence classes modulo $d'$. Since
$|A''| < k$ we have $d' > 1$. Let $S$ be the subset of elements $i \in \mathbb{Z}_{d'}$ such that
$A''$ contains the congruence class $i$ modulo $d'$. If $x$ is not a multiple of $d'$, and
$x = c'$ modulo $d'$, there exists an element $i \in S$ such that $i+c'$ modulo $d'$ is not in $S$,
otherwise $A''$ is the union of congruence classes modulo $gcd(c',d')$, contradicting the choice
of $d'$. This implies $A'$ does not contain any element from the congruence class $i$ modulo $d'$,
since otherwise $A'+x$ contains an element not in $A$. However, since $|A_2| \ge 2$, every congruence
class modulo $d$ that is contained in $A''$, and hence every congruence class modulo $d'$ that
is contained in $A''$, contains at least 2 elements in $A$. This implies there are two elements in
$A$ that are not in $A'$, contradicting the assumption that $|A'|=|A|-1$.

The remaining possibility is that $|A''| = k$. Since $|A| \le k-2$, this implies that there are at 
least two elements in the congruence class modulo $d$ that contains $A_2$, which are not contained 
in $A$. Let this class be $i$ modulo $d$. Suppose $x$ is not a multiple of $d$ and $x = c$ modulo
$d$. Since $A'+x \subset A$, at least two elements in the congruence class $i-c$ modulo $d$ are not
contained in $A'$. However, all elements in this class are contained in $A$, hence we again get
two elements in $A$ that are not contained in $A'$, a contradiction. Thus in this case $x$ must
be a multiple of $d$. In fact, using the argument in Case 1, we can show that in this case
$x$ must be $0, a$ or $-a$, but the weaker statement suffices.

This completes the proof of Lemma~\ref{inc}.
\hfill $\Box$

\begin{theorem}
\label{main}
Let $G$ be a $\mathbb{Z}_k$--weighted complete directed graph of order at least $k+2\Omega(k)$ for some 
$k \ge 2$. Then $G$ contains a  zero cycle.
\end{theorem}

We prove the theorem using the following lemma. 

\begin{lemma}
\label{one}
Let $G$ be a $\mathbb{Z}_k$--weighted complete directed graph of order at least $r+2\Omega(k)$, for some $1 
\le r < k$, $k \ge 2$, and let $u,v$ be any two distinct vertices in $G$.  Then either $G-\{u,v\}$ contains 
a zero cycle, or there exists a set $\{P_1,P_2,\ldots,P_r\}$ of $r$ $v$--$u$ paths in $G$, 
such that $|P_i| \ge 3$ and $w(P_i) \neq w(P_j)$, for all $1 \le i < j \le r$.
\end{lemma}

\noindent {\bf Proof:}[Theorem~\ref{main}]
Consider any $\mathbb{Z}_k$--weighted complete directed graph $G$ of order at least $k+2\Omega(k)$, for some 
$k \ge 2$. Let $x,y,z$ be any 3 vertices in $G$. If $w(xy) = w(xz)+w(zy)$ and $w(xz) = w(xy)+ w(yz)$ then adding 
the 2 equations gives $w(zy)+w(yz) = 0$, hence $z,y$ is a zero cycle. Without loss of generality, assume 
$w(xy) \neq w(xz) +w(zy)$. Applying Lemma~\ref{one} to the graph $G-\{z\}$ with $r = k-1$, either 
$G-\{x,y,z\}$ contains a zero cycle  or there exist $k-1$ $y$--$x$ paths of order at least 3 in $G-\{z\}$ 
with distinct weights. In the latter case, $G-\{z\}$ contains a $y$--$x$ path with weight $-w(xy)$ or a 
$y$--$x$ path with weight $-w(xz)-w(zy)$. In the first case, adding the edge $xy$ to the path gives a zero 
cycle, while in the second case, adding the edges $xz, zy$ gives a zero cycle. Thus in all cases $G$ contains 
a zero cycle.
\hfill $\Box$

Note that this shows that if Lemma~\ref{one} holds for some $k \ge 2$, then Theorem~\ref{main}
also holds for the same value of $k$.

\noindent {\bf Proof: }[Lemma~\ref{one}]
Suppose there exists a counterexample. Choose an example for which $k$ is minimum, and
subject to this condition $r$ is minimum.

If $r=1$, since $\Omega(k) \ge 1$, $|G| \ge 3$ and there exists a path of order at least 3 from $v$ to $u$.

Suppose $r=2$, which implies $|G| \ge 4$ and let $x,y$ be the vertices other than $u,v$ in $G$. If 
$w(xu) = w(xy) + w(yu)$ and $w(yu) = w(yx) + w(xu)$,  then adding the two equations gives $w(xy)+w(yx) = 0$. 
Thus $x,y$ is a zero cycle in $G-\{u,v\}$, a contradiction. Without loss of generality, assume $w(xu) \neq 
w(xy)+w(yu)$. Then the paths $v,x,u$ and $v,x,y,u$ have distinct weights, a contradiction.

Suppose that $r \ge 3$, and let $x$ be any vertex other than $u,v$.
By the minimality of $G$, either $G-\{u,x,v\}$ contains a zero cycle or there exist at least
$r-1$ $v$--$x$ paths of order at least 3 and distinct weights in $G-\{u\}$. We may assume the latter holds. 
If there are $r$ such paths with distinct weights, then appending the edge $xu$ to each of 
them gives $r$ $v$--$u$ paths in $G$ with distinct weights, a contradiction. We may 
assume $G$ contains $r-1$ $v$--$u$ paths of distinct weights. Let $A \subset \mathbb{Z}_k$ be 
the set of all elements $a \in \mathbb{Z}_k$, such that $G$ contains a $v$--$u$ path of weight $a$. 
Thus $|A| = r-1$. We show that $A$ must be a near arithmetic progression.
 
Let $(x,y)$ be any ordered pair of vertices in $G-\{u,v\}$. As argued in the case when $r=2$, we must
have either $w(xu) \neq w(xy)+w(yu)$ or $w(yu) \neq w(yx)+w(xu)$, otherwise the cycle $x,y$ has 0
weight. Again assume without loss of generality that $w(xu) \neq w(xy)+w(yu)$. By the minimality of $G$,  
either $G-\{u,x,y,v\}$ contains a zero cycle or there exists a set $\{P_1',P_2',\ldots,P_{r-2}'\}$ of 
$v$--$x$ paths in $G-\{u,y\}$, such that $|P_i'| \ge 3$ and $w(P_i') \neq w(P_j')$, for all $1 \le i < j 
\le r-2$. We may assume that the latter holds.

Let $Q_i$ be the $v$--$u$ path obtained by adding the edge $xu$ to $P_i'$ and let $S_i$ be the path obtained 
by adding edges $xy$ and $yu$ to $P_i'$.  Let $a = w(xy)+w(yu)-w(xu) \neq 0$, let $w(Q_i) = a_i$, and let 
$B = \{a_1,a_2,\ldots,a_{r-2}\}$. Note that $w(S_i) = a_i+a$.  We must have $B \cup (B+a) \subseteq A$ 
and therefore $B$ is a subset of $A$ with $|B| = |A|-1$ and $B+a \subset A$ for a nonzero $a \in \mathbb{Z}_k$.
Since $2 \le |A| \le k-2$, $A$ is a near arithmetic progression.  Let $x,y$ be any pair of vertices
in $G-\{u,v\}$. As argued previously, either $G-\{u,v,y\}$ contains a zero cycle or there exist $r-2$
$v$--$x$ paths in $G-\{u,y\}$ of order at least 3 with distinct weights. Adding the edge $xu$ to each such 
path we get a set of $r-2$ $v$--$u$ paths with distinct weights. Let $A' \subset A$ be the set of weights 
of these paths and let $e=w(xy)+w(yu)-w(xu)$. Replacing the edge $xu$ in these paths by the path $x,z,u$ gives 
a set of $v$--$u$ paths with distinct weights $A'+e$. Hence $A'+e \subset A$. Similarly, considering the 
set of $r-2$ $x$--$u$ paths in $G-\{v,y\}$, and $e = w(vy)+w(yx)-w(vx)$, we get a subset $A' \subset A$ with
$|A'| = |A|-1$ and $A'+e \subset A$. Therefore, by Lemma~\ref{inc}, one of the following properties
must hold.

\begin{enumerate}
\item
There exists a proper divisor $d > 1$ of $k$ such that for every ordered pair of vertices $x,y$ in $G-\{u,v\}$,
$w(xy)+w(yu)-w(xu)$ and $w(vy)+w(yx)-w(vx)$ is a multiple of $d$.
\item
There exists a nonzero constant $a \in \mathbb{Z}_k$ such that $gcd(a,k) = 1$ and for every ordered pair of 
vertices $x,y$ in $G-\{u,v\}$, $w(xy)+w(yu)-w(xu)$ and $w(vy)+w(yx)-w(vx)$ are contained in $\{0,a,-a\}$. 
\end{enumerate}

\noindent {\bf Case 1.} Suppose for all pair of vertices $x,y$ in $G-\{u,v\}$, $w(xy)+w(yu)-w(xu)$ and
$w(vy)+w(yx)-w(vx)$ are multiples of $d$ for some proper divisor $d > 1$ of $k$. 

We show that all $v$--$u$ paths in $G$ of order at least 3 have weights that are congruent modulo $d$.
First consider any two paths of length 2, $v,x,u$ and $v,y,u$. Then $w(vx)+w(xu)$ $\equiv$ $w(vy)+w(yx)+w(xu)$
$\equiv$ $w(vy)+w(yu)$  modulo $d$. Also, for any path $P = v,v_1,v_2,\ldots,v_l,u$ for $l \ge 2$, since
$w(v_iu) \equiv w(v_iv_{i+1})+w(v_{i+1}u)$ modulo $d$, it follows by induction that $w(P) \equiv
w(vv_1) + w(v_1u)$ modulo $d$. 

Since all $v$--$u$ paths in $G$ have weights that are congruent modulo $d$, there can be at most
$k/d$ such paths with distinct weights in $\mathbb{Z}_k$. Since there are $r-1$ such paths, we
must have $r \le k/d+1$.

Suppose $r \ge k/d$. Then $G$ is a complete directed graph of order $k/d+2\Omega(k) \ge k/d+2+2\Omega(k/d)$.
Thus $G-\{u,v\}$ is a complete directed graph of order at least $k/d+2\Omega(k/d)$. Define a weight function 
$w'$ on the edges of $G-\{u,v\}$ by $w'(xy) = i$ if $w(xy)+w(yu)-w(xu) = id$ (in $\mathbb{Z}_k$),
for some $0 \le i < k/d$. This gives a $\mathbb{Z}_{k/d}$-weighting of $G-\{u,v\}$ and by the minimality
of $k$, Lemma~\ref{one} and hence Theorem~\ref{main} holds for $k/d$. Therefore there exists
a cycle $C = v_1,v_2,\ldots,v_l$ in $G-\{u,v\}$ such that $w'(C) = 0$ in $\mathbb{Z}_{k/d}$. Since
$w'(C) = \sum_{i=1}^l (w(v_iv_{i+1}) + w(v_{i+1}u) - w(v_iu))/d$, where $v_{l+1} = v_1$, is a
multiple of $k/d$, $\sum_{i=1}^l w(v_iv_{i+1})$ is a multiple of $k$. Hence $C$ is also a zero
cycle in $\mathbb{Z}_k$ with weights $w$.

Suppose $r < k/d$. Again define the same $\mathbb{Z}_{k/d}$ weight function $w'$ for the edges $xy$
in $G-\{u,v\}$. Suppose for every $v$--$u$ path $P$ in $G$ of order at least 3, $w(P) \equiv a$ modulo $d$. 
Suppose for a vertex $x$ other than $u,v$ the weight of the path $v,x,u$ is $id+a$ for some
$0 \le i < k/d$. Define $w'(vx) = i$ and $w'(xu) = 0$. It is easy to show by induction on
the length that for any $v$--$u$ path $P$ in $G$ of order at least 3, $w(P) = w'(P)d+a$. Since 
$r < k/d$, applying Lemma~\ref{one} for $k/d$, either there exists a cycle $C$ in $G-\{u,v\}$ with
$w'(C) = 0$, or there exist $r$ $v$--$u$ paths in $G$ of order at least 3 with distinct weights 
in $\mathbb{Z}_{k/d}$. If the cycle $C$ exists, then $w(C) = 0$, and it is a zero cycle in $G-\{u,v\}$.
If the paths exist, then the same paths also have distinct weights in $\mathbb{Z}_k$, giving a contradiction
in either case. 

\noindent {\bf Case 2.} Suppose there exists a constant $a \in \mathbb{Z}_k$ such that $gcd(a,k) = 1$ and for 
all pairs of vertices $x,y$ in $G-\{u,v\}$, $w(xy)+w(yu)-w(xu)$ and $w(vy)+w(yx)-w(vx)$ are contained
in $\{0, a, -a\}$.

Let $w'(xy) = w(xy)+w(yu)-w(xu)$ and thus $w'(xy) \in \{0,a,-a\}$ for all edges $xy$ in $G-\{u,v\}$. 
If $C$ is a cycle in $G-\{u,v\}$, then $w'(C) = w(C)$, hence we may assume there is no cycle $C$ in 
$G-\{u,v\}$ such that $w'(C) = 0$. We will henceforth refer to $w'$ as the weight of an edge or
a subgraph, unless stated explicitly otherwise.

Suppose $G-\{u,v\}$ contains 3 vertices $x,y,z$ such that $w'(xy) = w'(yz) = c$ and $w'(xz) = -c$
for some $c \in \{a,-a\}$. We call such a triple of vertices a \emph{heavy} triple. If $r = 3$, consider 
the $v$--$u$ paths $P_1 = v,x,u$, $P_2 = v,x,y,u$ and $P_3 = v,x,y,z,u$. Then $w(P_2) = w(P_1) + c$ and 
$w(P_3) = w(P_1)+2c$, hence we get 3 $v$--$u$ paths of order at least 3 with distinct weights. 
If $r \ge 4$, by induction, either $G-\{u,v,x,y,z\}$ contains a cycle with $w(C) = 0$ or there exists 
$r-3$ $v$--$x$ paths in $G-\{u,y,z\}$ with distinct weights. Assume the latter holds. For any such path $P$, 
let $Q_1$ be the $v$--$u$ path obtained by adding the edge $xu$ to $P$, $Q_2$ is obtained by adding the edges 
$xy, yu$ to $P$, $Q_3$ is obtained by adding the edges $xz, zu$ to $P$, and $Q_4$ is obtained by adding the 
edges $xy, yz, zu$ to $P$. Then $w(Q_2) = w(Q_1)+c$, $w(Q_3) = w(Q_1)-c$ and $w(Q_4) = w(Q_1)+2c$. This implies 
that $A$ contains a subset $A'$ with $|A'| = r-3$ such that $A' \cup (A'-c) \cup (A'+c) \cup (A'+2c) \subseteq 
A$. Since $r-3 = |A'-c|$ $\le$ $|(A'-c) \cup A'|$ $\le$ $|(A'-c) \cup A' \cup (A'+c)|$ $\le |(A'-c) \cup A' 
\cup (A'+c) \cup (A'+2c)|$ $\le |A| = r-1$, one of the inequalities must be an equality. This implies there 
exists a subset $A'' \subset \mathbb{Z}_k$ such that $A''+c = A''$ for some $c$ with $gcd(c,k) = 1$. But this
implies $|A''| = k$, a contradiction. We may therefore assume $G$ does not contain a heavy triple of
vertices.

We now show that for some $c \in \{a,-a\}$, $G-\{u,v\}$ contains a Hamiltonian path $x_1,x_2,\ldots,x_l$
such that $w'(x_ix_{i+1}) = c$ for all $1 \le i < l$, where $l = |G-\{u,v\}|$. Since $|G| \ge r+2$, 
$l \ge r$,and we get $r$ $v$--$u$ paths in $G$ of order at least 3 and distinct weights. If $P_i$ is the 
$v$--$u$ path $v,x_1,\ldots,x_i,u$ for $1 \le i \le r$, it follows that $w(P_{i+1}) = w(P_i) + c$ and the 
paths have distinct weights. This gives a contradiction.

The edges $xy$ with $w'(xy) = 0$ must form a directed acyclic graph otherwise we get a zero cycle. If all 
the other edges have weight $c$ for some $c \in \{a,-a\}$, then it is easy it see that the graph has a 
Hamiltonian path with all edges having weight $c$. The vertices can be ordered $x_1,x_2,\ldots,x_l$ so that 
if $w'(x_ix_j) = 0$ then $i > j$. Then $w'(x_ix_{i+1}) = c$ for $1 \le i < l$ and this ordering gives the 
required Hamiltonian path.

Suppose $G-\{u,v\}$ has edges with weights $a$ as well as $-a$. We show that $G-\{u,v\}$ has a
very special structure, which implies it has a Hamiltonian path with all edges of weight $c$,
for some $c \in \{a,-a\}$. We claim that $G-\{u,v\}$ satisfies the following properties.
\begin{enumerate}
\item
The edges with weight 0 form a directed acyclic graph.
\item
There exists an edge $xy$ such that $w'(xy) = 0$ and for every vertex $z \not\in \{x,y\}$, either
$w'(zx) = 0$ or $w'(yz) = 0$. Note that both cannot hold since there is no cycle with 0
weight edges. We call such an edge the \emph{dominating} edge.
\item
For some $c \in \{a,-a\}$, $yx$ is the only edge in the graph with weight $-c$, all other edges have
weight $0$ or $c$.
\end{enumerate}

Note that these properties imply that $w'(xz) = w'(zy) = c$ for all vertices $z \not\in \{u,v,x,y\}$,
otherwise $x,z,y$ is a triangle with 0 weight.

We prove this by induction on $l$ for $l \ge 3$. Suppose $l = 3$. The edges of weight 0 must form a
directed acyclic graph. Suppose there are 3 edges with weight 0, say $xy$, $yz$ and $xz$. Then
the edges $yx$, $zy$ and $zx$ have nonzero weights and two of them must be the same. If $w'(zy)
\neq w'(yx)$ then the triangle $x,z,y$ has 0 weight, a contradiction. Therefore $w'(zy) = w'(yx)$
but $w'(zx) \neq w'(zy)$, which implies $(z,y,x)$ is a heavy triple, a contradiction.
 
Suppose there are two edges of 0 weight, say $xy, xz$. Then $w'(yz)$ and $w'(zy)$ are both nonzero 
and must be equal, say $w'(yz) = w'(zy) = c$ for some $c \in \{a,-a\}$. If either $w'(yx)$ or $w'(zx)$
is $-c$, then either $x,z,y$ or $x,y,z$ are triangles with 0 weight, respectively. A symmetrical argument
holds if $w'(yx) = w'(zx) = 0$. Suppose $w'(xy) = w'(yz) = 0$. Again, we must have $w'(xz) = w'(zx) =c$
for some $c \in \{a,-a\}$. If both $w'(yx) = w'(zy) = -c$, then $(z,y,x)$ is a heavy triple. Therefore
exactly one of $yx$ and $zy$ has weight $-c$ and the other has weight $c$. Then $G$ satisfies the required
properties with either $xy$ or $yz$ as the dominating edge.

If $xy$ is the only edge with weight 0, then we must have $w'(xz) = w'(zx)$ and $w'(yz) = w'(zy)$.
Also, $w'(zx) = w'(yz)$ otherwise the triangle $x,y,z$ has 0 weight. This implies we must have
$w'(yz) = w'(zx) = c$ and $w'(yx) = -c$, which implies $(y,z,x)$ is a heavy triple.

Finally, if there is no edge with 0 weight then $w'(xy) = w'(yx)$, $w'(xz) = w'(zx)$ and $w'(yz) = w'(zy)$.
Without loss of generality, we may assume $w'(xy) = w'(yz) = c$ and $w'(xz) = -c$ for some $c \in
\{a,-a\}$. This implies $(x,y,z)$ is a heavy triple.

Suppose $l \ge 4$. We claim that there exists a vertex $z$ in $G-\{u,v\}$ such that $G-\{u,v,z\}$
contains an edge of weight $a$ as well as an edge of weight $-a$. Let $x$ be a vertex in $G-\{u,v\}$
and suppose that $G-\{u,v,x\}$ has no edge of weight $-c$ for some $c \in \{a,-a\}$. We may assume
that either $xy$ or $yx$ is an edge with weight $-c$. Let $p$ and $z$ be two vertices other than
$x,y$ in $G-\{u,v\}$. At least one of the edges $yp$ or $py$ has weight $c$ and hence $G-\{u,v,z\}$
has an edge of weight $a$ as well as an edge of weight $-a$.

Applying induction, we may assume that $G-\{u,v,z\}$ has a dominating edge $xy$ of 0 weight, and
$yx$ is the only edge of weight $-c$ in $G-\{u,v,z\}$. 

Suppose all edges in the subgraph induced by $\{x,y,z\}$ have weights 0 or $-c$.
Suppose both $yz$ and $zx$ have weight $-c$. Let $p$ be any vertex other than $x,y,z$. If $w'(px) = 0$, 
then $w'(xp) = c$ and also $w'(py) = c$. This implies the cycle $x,p,y,z$ has 0 weight. If $w'(yp) = 0$, then 
$w'(py) = c$ and  $w'(xp) = c$, which again implies the cycle $x,p,y,z$ has 0 weight. Since $xy$ is a 
dominating edge, one of the two conditions must hold, and we get a contradiction. If both $yz$ and
$zx$ have weight 0, then $x,y,z$ is a triangle with weight 0. Suppose $w'(yz) = -c$ and $w'(zx) = 0$.
Then we have $w'(xz) = -c$. If $w'(px) = 0$, then if $w'(zp) = 0$ the triangle $p,y,z$ has 0
weight and if $w'(zp) = c$ the triangle $p,x,z$ has 0 weight. This implies $w'(zp) = -c$
and $(y,z,p)$ is a heavy triple. If $w'(yp) = 0$, then if $w'(zp) = 0$ the triangle $p,y,z$ has
0 weight and if $w'(zp) = c$, the cycle $p,y,x,z$ has 0 weight. Therefore $w'(zp) = -c$ and
hence $(x,z,y)$ is a heavy triple. This gives a contradiction. A symmetrical argument holds if $w'(yz) = 0$ 
and $w'(zx) = -c$. This implies $w'(zy) = -c$. If $w'(px) = 0$, then if $w'(pz) = 0$ the triangle $p,z,x$ 
has weight 0 and if $w'(pz) = c$, the cycle $p,z,y,x$ has weight 0. This implies $w'(pz) = -c$ and 
$(p,z,y)$ is a heavy triple. If $w'(yp) = 0$, then if $w'(pz) = 0$ the triangle $p,z,x$ has 0 weight and if 
$w'(pz) = c$ the triangle $p,z,y$ has 0 weight. Therefore $w'(pz) = -c$ and $(p,z,x)$ is a heavy triple.

We may assume that the subgraph induced by $\{x,y,z\}$ contains an edge of weight $a$ and
also an edge of weight $-a$. The argument for the base case implies there exists a dominating
edge in this subgraph, which can only be one of $xy$, $zx$ or $yz$.

Suppose $zx$ is the dominating edge. Then we have $w'(zx) = 0$, $w'(xz) = c$, and $w'(yz) = w'(zy)
= -c$. If $w'(xp) = 0$ then if $w'(zp) = 0$ the triangle $p,y,z$ has 0 weight and if $w'(zp) = c$,
the cycle $p,x,y,z$ has 0 weight. This implies $w'(zp) = -c$ and $(y,z,p)$ is a heavy triple.
If $w'(yp) = 0$ then if $w'(pz) = 0$ the cycle $p,z,x,y$ has 0 weight and if $w'(pz) = c$ the
triangle $p,z,y$ has 0 weight. Therefore $w'(pz) = -c$, which implies $(p,z,y)$ is a heavy triple.

Suppose $yz$ is the dominating edge. Then we have $w'(yz) = 0$, $w'(zy) = c$, $w'(xz) = w'(zx) = -c$. 
If $w'(px) = 0$ then if $w'(zp) = 0$ the cycle $p,x,y,z$ has 0 weight and if $w'(zp) = c$ the triangle
$p,x,z$ has 0 weight. Therefore $w'(zp) = -c$ and $(x,z,p)$ is a heavy triple. If $w'(yp) = 0$
then if $w'(pz) = c$ the cycle $z,x,y,p$ has 0 weight and if $w'(pz) = -c$ the triangle $z,y,p$
has weight 0. Therefore $w'(pz) = 0$. Similarly, if $w'(zp) = c$, the cycle $z,p,y,x$ has 0 weight
and if $w'(zp) = -c$ the triangle $z,p,y$ has 0 weight. Therefore $w'(zp) = 0$ and the cycle $p,z$
has 0 weight, a contradiction.

Therefore the dominating edge must be $xy$ itself. This implies exactly one of the edges $zx$, $yz$
has 0 weight and all others, except $yx$, have weight $c$ in the subgraph induced by $\{x,y,z\}$. Suppose 
$w'(zx) = 0$. We claim that for any vertex $p \not\in \{x,y,z\}$ neither the edge $zp$ nor the edge $pz$
can have weight $-c$. If $w'(px) = 0$ then if $w'(zp) = -c$  then $(z,y,p)$ is a heavy triple and
if $w'(pz) = -c$ then $(p,y,z)$ is a heavy triple. If $w'(yp) = 0$ then if $w'(zp) = -c$ the
cycle $z,p,y,x$ has 0 weight and if $w'(pz) = -c$ the triangle $z,y,p$ has 0 weight. Therefore $yx$
is the only edge with weight $-c$ in  $G-\{u,v\}$ and $G-\{u,v\}$ satisfies all the required properties.
Suppose $w'(yz) = 0$. If $w'(px) = 0$ then if $w'(zp) = -c$ the triangle $p,x,z$ has 0 weight and
if $w'(pz) = -c$ the cycle $p,z,y,x$ has 0 weight. If $w'(yp) = 0$, then if $w'(zp) = -c$ the
cycle $p,y,x,z$ has weight 0 and if $w'(pz) = -c$ the triangle $p,z,y$ has weight 0. Therefore
$yx$ is the only edge with weight $-c$, and $G-\{u,v\}$ satisfies all the required properties.

Now it is easy to show that $G-\{u,v\}$ has a Hamiltonian path with all edges of weight $c$.
Let $x_1,x_2,\ldots,x_l$ be an ordering of the vertices such that if $w'(x_ix_j) = 0$ then
$i > j$. Since $xy$ is a dominating edge, we must have $y = x_i$ and $x = x_{i+1}$ for
some $1 \le i < l$.  Suppose $i > 1$ and $i+1 < l$.  Then $x_1,x_2,\ldots,x_{i-1},x_{i+1},\ldots,x_l,x_i$
is a Hamiltonian path with all edges of weight $c$. If $i = 1$, since $l \ge 3$, the path 
$x_2,x_3,\ldots,x_l,x_1$ is a Hamiltonian path with all edges of weight $c$. Similarly, if $i+1 = l$,
the path $x_l,x_1,x_2,\ldots,x_{l-1}$ is the required Hamiltonian path.

This completes the proof of Lemma~\ref{one} and hence the proof of Theorem~\ref{main}.
\hfill $\Box$

\section{Undirected Graphs}

In this section, we consider undirected graphs in which weights are assigned to
vertices as well as edges. The weight of a subgraph is the sum of the weights of vertices
and edges in the subgraph. The proof of Theorem~\ref{main} can be modified slightly to show that every 
$\mathbb{Z}_k$-complete weighted directed graph of order at least $k+1+2\Omega(k)$ contains a directed zero cycle 
of length at least 3. Lemma~\ref{one} also needs to be modified slightly to consider graphs of
order at least $r+1+2\Omega(k)$ to ensure the cycle has length at least 3. If $r = 2$ and $k > 2$, 
the statement fails for a complete graph of order 4, if $w(xy) = w(yx) = 0$,
$w(vx) = w(vy) = w(xu) = w(yu)= 1$. However, the induction step is the same with the hypothesis
that the zero cycle in $G-\{u,v\}$ has length at least 3. The following statement therefore
follows from the modified Theorem~\ref{main}, after replacing each undirected edge by two 
oppositely directed edges with the same weight, and adding the weight of a vertex to all edges
directed out of the vertex. Since the zero cycle has length at least 3, it is also a cycle in
the undirected graph.

\begin{corollary}
Let $G$ be a complete undirected graph of order $k+1+2\Omega(k)$ for some $k \ge 2$ and suppose
every vertex and edge in $G$ is assigned a weight in $\mathbb{Z}_k$. Then $G$ contains a cycle with
weight 0.
\end{corollary}

More significantly, in the undirected case, a similar result holds for all graphs with
sufficiently large minimum degree, rather than just complete graphs. The bound that
we prove is weaker though.

\begin{theorem}
\label{undirected}
Let $\Gamma$ be any nontrivial finite abelian group and $G$ an undirected 
$\Gamma$-weighted graph with minimum degree at least $2|\Gamma|-1$.
Then there exists a zero cycle in $G$.
\end{theorem}

Let $|\Gamma| = k \ge 2$ and denote by $w(x)$ the weight of a vertex,
edge or subgraph of $G$. The proof technique used is almost exactly the
same as used in~\cite{D} to prove a completely different result. Let $N(x)$
denote the set of vertices in $G$ that are adjacent to the vertex $x$ in $G$.

We consider ordered pairs of the form $(G,K)$ where $K$ is a
proper complete subgraph of the graph $G$, possibly empty.

\begin{definition}
The ordered pair $(G,K)$ is said to contain a configuration of type 
\textbf{A} if there exists a vertex $x \in V(G) \setminus V(K)$ such that 
$|N(x) \cap V(K)| \ge 2k-1$.
\end{definition}

\begin{definition}
The ordered pair $(G,K)$ is said to contain a configuration of type 
\textbf{B} if there exist two vertices $x, y$ in $V(G) \setminus V(K)$ such
that $|N(x) \cap V(K)| \ge 2k-2$, $|N(y) \cap V(K)| \ge 2k-2$ and there
exists an $x$--$y$ path $P$ in $G-V(K)$.
\end{definition}

\begin{definition}
The ordered pair $(G,K)$ is said to contain a configuration of type
\textbf{C} and rank $r$, $2 \le r \le k$, if there exist two vertices $x,y$ in $V(G) \setminus
V(K)$ such that $|N(x) \cap V(K)| \ge 2(k-r)+1$, $|N(y) \cap V(K)|
\ge 2(k-r)+1$ and there exists a set $\{P_1, P_2, \ldots, P_r\}$ of $r$
$x$--$y$ paths in $G-V(K)$ having distinct weights.
\end{definition}

\begin{definition}
The ordered pair $(G,K)$ is said to contain a configuration of type \textbf{D} and \emph{rank} $r$, 
$2 \le r < k$, if $V(G)\setminus V(K)$ contains three vertices $x,z,y$ satisfying the following properties.
\begin{enumerate}
\item
$xz$ and $yz$ are edges in $G-V(K)$.
\item
There exist vertices $x', y', z' \in V(G) \setminus V(K)$ and vertex disjoint paths $P_x, P_y, P_z$ in 
$G-V(K)$ such that $|N(x') \cap V(K)| \ge 2(k-r)-1$, $|N(y') \cap V(K)| \ge 2(k-r)$, 
$|N(z') \cap V(K)| \ge 2(k-r)$, and for all $u \in \{x,y,z\}$, $P_u$ is a $u$--$u'$ path. Note that the 
vertex $u'$ may be the same as the vertex $u$, in which case the path $P_u$ is trivial, for $u \in \{x,y,z\}$.
\item
There is a set $\{P_1,P_2,\ldots,P_r\}$ of $r$ $x$--$y$ paths in $G-V(K)$ having distinct weights 
such that $P_i$ is internally vertex disjoint from $P_x$, $P_y$ and $P_z$ for $1 \le i \le r$.
\end{enumerate}
\end{definition}

\begin{lemma}
\label{reduction}
Let $G$ be a $\Gamma$-weighted graph and $K$ a proper complete subgraph of $G$. If every vertex in 
$V(G)\setminus V(K)$ has degree at least $2k-1$ in $G$, then either $G-V(K)$ contains a zero cycle, 
or the ordered pair $(G,K)$ contains a configuration of one of the types \textbf{A}, \textbf{B}, \textbf{C}, 
or \textbf{D}. 
\end{lemma}

\noindent {\bf Proof:} Let $(G,K)$ be a counterexample that minimizes $|V(G)|+|V(G)\setminus V(K)|$. If  
$|V(G) \setminus V(K)| = 1$, the only vertex $x \in V(G) \setminus V(K)$ has degree at least $2k-1$, hence 
$|N(x) \cap V(K)| \ge 2k-1$. This implies $(G,K)$ contains a configuration of type \textbf{A}, a contradiction.

Suppose $|V(G) \setminus V(K)| > 1$. We consider two cases, one of which is straightforward.

\noindent {\bf Case 1.}
Suppose there exists a vertex $v \in V(G) \setminus V(K)$ that is adjacent to all vertices in $V(K)$. Let $K'$ 
be the complete subgraph of $G$ induced by $V(K) \cup \{v\}$. Then the ordered pair $(G,K')$ satisfies the 
hypothesis of Lemma \ref{reduction}, and by the minimality of $(G,K)$,  either $G-V(K')$ contains a zero cycle 
or $(G,K')$ contains a configuration of one of the four types. Since $G-V(K')$ is a subgraph of $G-V(K)$, we 
may assume that the latter holds. Now we show that in each case, the configuration in $(G,K')$ can be modified 
to either find a zero cycle in $G-V(K)$, or get a configuration of one of the four types in $(G,K)$, contradicting 
the fact that $(G,K)$ is a counterexample.

\noindent {\bf Case 1.1}
Suppose $(G,K')$ contains a configuration of type \textbf{A}. Let $x$ be a vertex in $V(G) \setminus V(K')$ 
such that $|N(x) \cap V(K')| \ge 2k-1$. If $x$ is not adjacent to $v$, then $|N(x) \cap V(K)| \ge 2k-1$ and 
$(G,K)$ contains a configuration of type \textbf{A}. If $x$ is adjacent to $v$, then $|N(x) \cap V(K)| \ge 2k-2$ 
and since $v$ is adjacent to every vertex in $V(K)$, $|N(v) \cap V(K)| = |V(K)| \ge 2k-2$. The edge $vx$ 
implies that $(G,K)$ contains a configuration of type \textbf{B}.

\noindent {\bf Case 1.2}
Suppose $(G,K')$ contains a configuration of type \textbf{B}. Let $x,y$ be vertices in $V(G) \setminus V(K')$ 
such that $|N(x) \cap V(K')| \ge 2k-2$, $|N(y) \cap V(K')| \ge 2k-2$ and there is an $x$--$y$ path $P$ in 
$G-V(K')$.

If $v$ is not adjacent to both the vertices $x,y$, then $x,y$ satisfy the same properties with $K'$ replaced by 
$K$, and $(G,K)$ contains the same configuration of type \textbf{B}. 

If $v$ is adjacent to $x$ but not adjacent to $y$, then $|N(y) \cap V(K)| \ge 2k-2$ and hence 
$|N(v) \cap V(K)| \ge 2k-2$. Also, $P \cup vx$ is a $v$--$y$ path in $G-V(K)$. Hence $(G,K)$ contains 
a configuration of type \textbf{B}.  A symmetrical argument holds if $v$ is adjacent to $y$ but not adjacent 
to $x$.

Suppose $v$ is adjacent to both $x$ and $y$. Then $|N(x) \cap V(K)| \ge 2k-3$, $|N(y) \cap V(K)| \ge 2k-3$ and 
hence $|N(v) \cap V(K)| \ge 2k-3$. If $w(x)+w(vx)+w(v)+w(vy)+w(y) \neq w(P)$ then $P$ and $vx \cup vy$ are two 
$x$--$y$ paths in $G-V(K)$ with distinct weights. This implies $(G,K)$ contains a configuration of type
\textbf{C} and rank two. Similarly, if either $w(v)+w(vx)+w(P) \neq w(v)+w(vy)+w(y)$ or $w(P)+w(vy)+w(v) 
\neq w(x)+w(vx)+w(v)$ then $(G,K)$ contains a configuration of type \textbf{C} and rank two. If none of these 
inequalities holds, then the three equations imply that the cycle $P \cup vx \cup vy$ is a zero cycle. 

Hence, either $G-V(K)$ contains a configuration of type \textbf{C} and rank two, or $G-V(K)$ contains a zero
cycle, a contradiction.

\noindent {\bf Case 1.3}
Suppose $(G,K')$ contains a configuration of type \textbf{C} and rank $r$, for some $2 \le r \le k$. 
Let $x,y$ be vertices in $V(G) \setminus V(K')$ such that $|N(x) \cap V(K')| \ge 2(k-r)+1$, $|N(y) \cap V(K')| 
\ge 2(k-r)+1$ and let $\{P_1,P_2,\ldots,P_r\}$ be the set of $r$ $x$--$y$ paths in $G-V(K')$ having distinct 
weights.

If $v$ is not adjacent to both $x$ and $y$, then $(G,K)$ also contains the same configuration of type 
\textbf{C} and rank $r$.

If $v$ is adjacent to $x$ but not to $y$, then $|N(y) \cap V(K)| \ge 2(k-r)+1$ and hence $|N(v) \cap V(K)| 
\ge 2(k-r)+1$. The paths $P_i \cup vx$ for $1 \le i \le r$ are $v$--$y$ paths in $G-V(K)$ having distinct 
weights. Hence$(G,K)$ contains a configuration of type \textbf{C} and rank $r$. A symmetrical argument holds 
if $v$ is adjacent to $y$ but not adjacent to $x$.

Suppose $v$ is adjacent to both $x$ and $y$. If $r = k$, then the cycles 
$vx \cup vy \cup P_i$ for $1 \le i \le k$ have distinct weights, and since $|\Gamma| = k$, one of these has 
weight 0. Therefore $G-V(K)$ contains a zero cycle. Suppose $r < k$. Since $|N(x) \cap V(K)| \ge 2(k-r)$, 
$|N(y) \cap V(K)| \ge 2(k-r)$ we have $|N(v) \cap V(K)| \ge 2(k-r)$. Relabeling the vertex $v$ as $z$, 
choosing the vertices $x', y', z'$ to be the vertices $x,y,z$ respectively and the paths 
$P_x, P_y, P_z$ to be trivial, we get a configuration of type \textbf{D} and rank $r$ in $(G,K)$.

\noindent {\bf Case 1.4}
Finally, suppose $(G,K')$ contains a configuration of type \textbf{D} and rank $r$, for some $2 \le r < k$.
Let $x, y, z$ be the three vertices in $V(G) \setminus V(K')$ that satisfy the properties defined in
configuration \textbf{D}, such that $xz, yz$ are edges in $G$. Let $x', y', z'$ be the vertices in $V(G) 
\setminus V(K')$ such that $|N(x') \cap V(K')| \ge 2(k-r)-1$, $|N(y') \cap V(K')| \ge 2(k-r)$ and 
$|N(z') \cap V(K')| \ge 2(k-r)$ and let $P_x, P_y, P_z$ be the vertex disjoint $x$--$x'$, $y$--$y'$ and 
$z$--$z'$ paths in $G-V(K')$, respectively.  Let $\{P_1,P_2, \ldots, P_r\}$ be the set of $r$ $x$--$y$ paths 
in $G-V(K')$ that are internally vertex disjoint from the paths $P_x, P_y, P_z$ and have distinct weights.

If $v$ is not adjacent to any of the vertices $x',y',z'$, it is clear that $(G,K)$ contains the same 
configuration of type \textbf{D} and rank $r$.  If $v$ is adjacent to exactly one of the vertices $x', y', z'$, 
then $|N(v) \cap V(K)| \ge 2(k-r)$. If $v$ is adjacent only to $u'$ for some $u \in \{x,y,z\}$, replace the 
vertex $u'$ by the vertex $v$ and the path $P_u$ by the path $P_u \cup u'v$. This gives a configuration of 
type \textbf{D} and rank $r$ in $(G,K)$. 

Suppose $v$ is adjacent to $x'$ and $y'$ but not to $z'$. Then $|N(z') \cap V(K)| \ge 2(k-r)$, 
$|N(y') \cap V(K)| \ge 2(k-r)-1$ and hence $|N(v) \cap V(K)| \ge 2(k-r)$. Replace the vertex $x'$ by 
$v$ and the path $P_x$ by the path $P_x \cup vx'$. Now interchange the labels of the vertices $x,y$, label $v$ 
as $y'$ and $y'$ as $x'$, to get a configuration of type \textbf{D} and rank $r$ in $(G,K)$. 

Suppose $v$ is adjacent to $x'$ and $z'$ and $v$ may or may not be adjacent to $y'$. Then $|N(z') \cap V(K)| 
\ge 2(k-r)-1$, $|N(y') \cap V(K)| \ge 2(k-r)-1$ and hence $|N(v) \cap V(K)| \ge 2(k-r)-1$. Let 
$Q_i =  vx' \cup P_x \cup P_i \cup P_y$ and $Q_i' = vz' \cup P_z \cup xz \cup P_i \cup P_y$, for $1 \le i \le r$, 
be $v$--$y'$ paths in $G-V(K)$. If among the $2r$ paths $\{Q_1, Q_1',\ldots,Q_r,Q_r'\}$, there are $r+1$ paths 
of distinct weights, then $(G,K)$ contains a configuration of type \textbf{C} and rank $r+1$, with $v$ and $y'$ 
as the two required vertices satisfying the properties defined for configuration \textbf{C}. Similarly, 
let $S_i = vz' \cup vx' \cup P_x \cup P_i \cup P_y$ and $S_i' = P_z \cup 
xz \cup P_i \cup P_y$, for $1 \le i \le r$ be $z'$--$y'$ paths in $G-V(K)$. If among the $2r$ paths 
$\{S_1, S_1', \ldots, S_r,S_r'\}$ there are $r+1$ paths of distinct weights, then $(G,K)$ contains a 
configuration of type \textbf{C} and rank $r+1$, with $z'$ and $y'$ as the required vertices.

Suppose both sets of paths $\{Q_1,Q_1',\ldots,Q_r,Q_r'\}$ and $\{S_1,S_1',\ldots,S_r,S_r'\}$ contain only 
$r$ paths of distinct weights. Note that $w(Q_i) = w(P_i) + c_1$ for some $c_1 \in \Gamma$ and all 
$1 \le i \le r$, which implies that $\{Q_1,\ldots,Q_r\}$ have distinct weights. Similarly, 
$w(Q_i') = w(P_i)+c_2$ for some $c_2 \in \Gamma$.  Also $w(S_i) = w(Q_i)+w(z')+w(vz')$ and $w(S_i') = 
w(Q_i')-w(v)-w(vz')$.  Suppose  $w(P_i) = a$ for some $1 \le i \le r$ and $a \in \Gamma$. Then 
$w(Q_i') = a+c_2$ and there exists an index $j$ such that $w(Q_j) = a+c_2$. Therefore $w(S_j) = a+c_2+w(z')
+w(vz')$. Hence, there is an index $m$ such that $w(S_m') = a+c_2+w(z')+w(vz')$, which 
implies $w(Q_m') = a+c_2+w(z')+2w(vz')+w(v)$ and hence $w(P_m) = a + c$ where $c = w(z')+2w(vz')+ w(v)$. 
Since this holds for any path $P_i$, it implies that the set $\{P_1,\ldots,P_r\}$ contains a path of weight 
$a + ic$, for all $i \ge 0$.

Let $T_i$ be the $v$--$z'$ path $vx' \cup P_x \cup P_i \cup yz \cup P_z$ for $1 \le i \le r$. If none of the 
paths $T_i$ for $1 \le i \le r$, has weight $w(v) + w(vz') + w(z')$ then $\{T_1,T_2,\ldots,T_r\} \cup \{vz'\}$ 
is a set of $r+1$ $v$--$z'$ paths in $G-V(K)$ having distinct weights. This implies $G-V(K)$ contains a 
configuration of type \textbf{C} and rank $r+1$. Suppose, without loss of generality, that $w(T_1) = w(v) + 
w(vz') + w(z')$. By the previous discussion, there exists a path $P_i$ for some $1 \le i \le r$, such that 
$w(P_i) = w(P_1) - c$ and hence $w(T_i) = w(T_1) - c = -w(vz')$. This implies that the cycle $T_i \cup vz'$ 
is a zero cycle. Thus either $(G,K)$ contains a configuration of type \textbf{C} and rank $r+1$ or $G-V(K)$ 
contains a zero cycle.

If $v$ is adjacent to $z'$ and $y'$ but not to $x'$, then $|N(z') \cap V(K)| \ge 2(k-r)-1$, 
$|N(x') \cap V(K)| \ge 2(k-r)-1$ and hence $|N(v) \cap V(K)| \ge 2(k-r)-1$. Now, we can use the same argument 
as before, by interchanging the vertices $x,y$ and $x',y'$.

\noindent {\bf Case 2.}
Suppose there is no vertex $v$ in $V(G) \setminus V(K)$ that is adjacent to all vertices in $V(K)$.
Let $u$ be any vertex in $V(K)$. For every vertex $v \in N(u) \setminus V(K)$, let $f(v)$ denote
any vertex in $V(K)$ that is not adjacent to $v$. Let $G'$ be the graph obtained from $G-\{u\}$ by 
adding edges $vf(v)$ for all vertices $v \in N(u)\setminus V(K)$. Let $K' = K-\{u\}$.

Now, $|V(G') \setminus V(K')| = |V(G) \setminus V(K)|$ but $|V(K')| < |V(K)|$, hence by the minimality of 
$(G,K)$, either $G'-V(K')$ contains a zero cycle, or $(G',K')$ contains one of the four types of 
configurations. Since $G'-V(K') = G-V(K)$ and $|N(v) \cap V(K')|$ in $G'$ equals $|N(v) \cap V(K)|$ in $G$
for every vertex $v \in V(G)\setminus V(K)$, it follows from the definitions of the configurations that 
either $G-V(K)$ contains a zero cycle or $(G,K)$ contains the same configuration as $(G',K')$.

This completes the proof of Lemma \ref{reduction}.
\hfill $\Box$

\noindent {\bf Proof:}[Theorem~\ref{undirected}]
The proof of Theorem \ref{undirected} follows immediately from Lemma \ref{reduction}. If $G$ is a graph 
with minimum degree at least $2k-1$, then the ordered pair $(G,\emptyset)$ satisfies the hypothesis of 
Lemma \ref{reduction} and hence either $G$ contains a zero cycle or $(G,\emptyset)$ contains a configuration 
of one of the four types. However, since $K$ is empty, the latter is not possible, and the theorem follows.

\hfill $\Box$

\section{Remarks}

In the proof of Lemma~\ref{one}, in Case 2, we used the fact that $G-\{u,v\}$ does not
contain a heavy triple. However, it is possible that the conclusion is true even without this
property.

\noindent{\bf Question 1:} Let $G$ be a complete directed graph in which each edge is
assigned a weight in $\{0,1,-1\}$. Is it true that either $G$ contains a cycle with 0
weight or there exists a Hamiltonian path in $G$ with all edges having the same nonzero weight. 
Note that addition here is ordinary integer addition.

It is easy to construct a $\mathbb{Z}_k$-weighted complete directed graph of order $k$ that does
not contain a zero cycle~\cite{MS}. A similar construction gives a lower bound on the 
minimum degree of an undirected $\mathbb{Z}_k$-weighted graph that does not contain a zero cycle. 
Let $G$ be a graph obtained from a nontrivial tree by adding a complete graph of order $k-1$ 
and joining every vertex in the complete graph to every vertex in the tree. The minimum degree 
of this graph is $k$. Assign weight $1$ to the vertices in the complete subgraph and weight 0 
to all other vertices and all edges in the graph. Again this graph does not contain a zero cycle.

\noindent {\bf Question 2:} Is it true that every $\mathbb{Z}_k$-weighted undirected graph with
minimum degree $k+1$ contains a zero cycle? 

The extremal example in the undirected case also gives a lower bound $kn- k(k+1)/2$ 
on the number of edges. Since any graph with $n \ge k+1$ vertices and more than 
$kn-k(k+1)/2$ edges contains a subgraph with minimum degree $k+1$, a weaker form of Question 2 
is whether every such $\mathbb{Z}_k$-weighted graph contains a zero cycle. If so, the bound on the
number of edges is optimal for $n \ge k+1$.

Thomassen~\cite{T1} showed that there are directed graphs with arbitrarily large in and out degrees
that do not contain an even cycle. However, when restricted to strongly connected directed graphs,
he showed that every such $\mathbb{Z}_2$-weighted directed graph with minimum in and out degree 3 
contains a zero cycle~\cite{T2}. This leads to the following question.

\noindent {\bf Question 3.} Is it true that there exists a function $f(k)$ such that every
$\mathbb{Z}_k$-weighted strongly connected directed graph with minimum in and out degree at least $f(k)$
contains a zero cycle. If this is not true, does assuming large enough connectivity ensure 
the existence of such a cycle?

\end{document}